\documentclass[]{article}
\usepackage{graphicx}
\usepackage{graphics}
\usepackage{amssymb,amsmath,amsthm,epsfig}
\usepackage{color,soul}
\newcommand{\J}{P^{(\alpha,\alpha)}}
\newcommand{\wJ}{\widetilde P^{(\alpha,\alpha)}}

\newcommand\dint{\displaystyle\int}
\newcommand\ds{\displaystyle}
\newtheorem{theorem}{Theorem}
\newtheorem{lemma}{Lemma}

\newtheorem{proposition}{Proposition}

\newtheorem{remark}{Remark}
\newcounter{reh}
\newcounter{rek}
\begin{document}
\begin{center}
{\large {\bf Weighted Finite Laplace Transform Operator: Spectral Analysis  and Quality of Approximation by its Eigenfunctions.}}\\
\vskip 1cm NourElHouda Bourguiba$^a$   and Abderrazek Karoui$^a$ {\footnote{
Corresponding author: Abderrazek Karoui, Email: abderrazek.karoui@fsb.rnu.tn
This work was supported by the DGRST research Grant UR13ES47.}}
\end{center}
\vskip 0.5cm {\small

\noindent $^a$ University of Carthage,
Department of Mathematics, Faculty of Sciences of Bizerte, Tunisia.
}\\

\noindent{\bf Abstract}--- For  two  real numbers $c>0, \alpha> -1,$ we   study some spectral properties of the weighted finite bilateral Laplace transform operator, defined over the space $E=L^2(I,\omega_{\alpha}),$ $I=[-1,1],$ $\omega_{\alpha}(x)=(1-x^2)^{\alpha},$  by ${\displaystyle 
\mathcal L_c^{\alpha} f(x)= \int_I e^{cxy} f(y) \omega_{\alpha}(y)\, dy}.$ In particular, we use a technique based on the Min-Max theorem to  prove that the sequence of the eigenvalues of this  operator has a super-exponential decay rate to zero.  Moreover, we give a lower bound with a  magnitude of
order $e^c,$ for the largest eigenvalue of the operator $\mathcal L_c^{\alpha}.$    Also, we give some local estimates and bounds of the eigenfunctions $\varphi_{n,c}^{\alpha}$ of $\mathcal L_c^{\alpha}.$ Moreover, we show that these eigenfunctions are good candidates for the spectral approximation of a function that can be written as a weighted finite Laplace transform of an other $L^2(I,\omega_{\alpha})-$function.
 Finally, we give some numerical examples that  illustrate the different results of this work. In particular, we provide an example that illustrate the Laplace based spectral method, for the inversion of the finite Laplace transform.\\



\section{Introduction} It is well known that the finite Laplace transform has been used as a tool for solving various problems
from different scientific area such as Engineering, Physics, Mathematics, to cite but a few, see for example  
\cite{Dakto, Rutily, Debnath, Bertero, Ledermann},  some of these applications. To the best of our knowledge, H. S. Dunn was the first to define in \cite{Dunn},  finite Laplace transform over a finite interval of the time domain. He has also pointed out that the finite Laplace transform yields identical results as the usual Laplace transform. More importantly the finite version of this transform has the advantage to be more general, it can be applied to functions which are not of exponential order and it is applicable to solve boundary values problems.\\

We should  mention that in \cite{Debnath}, the authors have given the basic properties as well as an iversion formula for the finite Laplace transform over the finite interval $[0,T]$. Recently, in \cite{Ledermann},  the  authors have used some decay properties of the eigenvalues of the finite Laplace transform $\mathcal L_{a,b}$, over the interval $[a,b],$ as well as a growth condition on the eigenvalues  of the differential operator, commuting with $\mathcal L_{a,b}$ and they have given a lower bound for a norm of $\|\mathcal L_{a,b}\|.$ Note that the commutativity of the finite Laplace transform with the corresponding differential operator of the Sturm-Liouville type, has been first studied in \cite{Grunbaum}.

 In this work, we study  some spectral properties and computational techniques of the eigenvalues and eigenfunctions  of a fairly  more general setting of the weighted finite Laplace transform, defined by 
\begin{equation}\label{Eq1.1}
\mathcal L_c^{\alpha} f(x)= \int_I e^{cxy} f(y) \omega_{\alpha}(y)\, dy,\quad \omega_{\alpha}(y)=(1-y^2)^{\alpha},
\end{equation}
where, $c>0$ and $\alpha >-1$ are two real numbers and $f\in L^2(I,\omega_{\alpha}),$ $I=[-1,1].$ Note that for the special case where $\alpha=0,$ the operator $\mathcal L_c^{0},$ is reduced to the usual finite Laplace transform. Also, by a simple dilation and translation arguments, the results of this paper still hold for the general case of weighted finite Laplace transform, given by 
$$ L_c^{\alpha}(g)(x)= \int_a^b e^{cxy} g(y) \left((y-a)(b-y)\right)^{\alpha} \, dy,\quad a,b\in\mathbb R.$$
 Note that the previous  weighted finite Laplace  transform can be considered as the compactly supported weight function case of a more general
weighted Laplace transform operator, considered in   \cite{Bertero},  where the authors have shown that the weighted version of the Laplace transform
has the advantage to avoid edge effects in the numerical inversion of the Laplace transform.

In this work, we are interested in the properties of the eigenfunctions $\varphi_{n,c}^{\alpha}$ and the associated eigenvalues  $\nu_{n,\alpha}(c),$ $n\geq 0,$  of the operator $\mathcal L_c^{\alpha},$ which are solutions of the integral equation
\begin{equation}\label{Eq1.2}
\mathcal L_c^{\alpha} \varphi_{n,c}^{\alpha}(x)= \int_I e^{cxy} \varphi_{n,c}^{\alpha}(y) \omega_{\alpha}(y)\, dy = \nu_{n,\alpha}(c) \varphi_{n,c}^{\alpha}(x),\quad x\in I.
\end{equation}
It can be easily checked that the eigenvalues $\nu_{n,\alpha}(c)$ are simple. Moreover, 
in the sequel, we assume that these eigenvalues  are arranged in the decreasing order, so that 
$$ \nu_{0,\alpha}(c) > \nu_{1,\alpha}(c) >\cdots > \nu_{n,\alpha}(c)> \cdots.$$ Since $\mathcal L_c^{\alpha}$ is a self-adjoint operator, then  the eigenfunctions $\varphi_{n,c}^{\alpha}$ are orthogonal in $L^2(I,\omega_{\alpha}).$ Also, we assume that 
these eigenfunctions  are normalized so that 
\begin{equation}\label{normalization}
\int_I  \varphi_{n,c}^{\alpha}(y) \varphi_{m,c}^{\alpha}(y) \omega_{\alpha}(y)\, dy = \delta_{n,m}.
\end{equation}
We should mention that in \cite{Karoui-Souabni, Wang2}, a spectral analysis of the similar weighted finite Fourier transform has been studied. This last operator is given by 
\begin{equation}\label{Eq1.3}
\mathcal F_c^{\alpha} f(x)= \int_I e^{icxy} f(y) \omega_{\alpha}(y)\, dy,\quad \omega_{\alpha}(y)=(1-y^2)^{\alpha},\quad
c>0,\,\,\, \alpha>-1.
\end{equation}

Note that the  change of  the trigonometric exponential kernel of the operator $\mathcal F_c^{\alpha}$ to the real exponential kernel for the operator $\mathcal L_c^{\alpha},$ give rise to different spectral properties. For example, the real exponential kernel is a positive definite kernel and consequently, $\nu_{n,\alpha}(c) > 0.$ This is not the case for the eigenvalues $\mu_{n,\alpha}(c)$ of $\mathcal F_c^{(\alpha)},$ that are either purely imaginary complex numbers or real numbers and satisfy the equality ${\displaystyle 
\mu_{n,\alpha}(c)= i^n |\mu_{n,\alpha}(c)|}.$ Another main difference is that unlike the $\nu_{n,\alpha}(c),$ the eigenvalues $\mu_{n,\alpha}(c)$ for $\alpha\geq 0,$  are related to an energy maximization problem and as a consequence, we have for $
{\displaystyle \frac{2c}{\pi} |\mu_{0,\alpha}(c)|^2 < 1},$ for more details, see \cite{Karoui-Souabni}. This is not the case for  the largest eigenvalues 
$\nu_{0,\alpha}(c)$ which grows like $e^c.$ Nonetheless, we show that as for the case of the $\mu_{n,\alpha}(c),$ the sequence $\nu_{n,\alpha}(c)$ has a super-exponential decay rate. This is one of the main results of this work. The proof of this last result is completely different from the one used in \cite{Karoui-Souabni} and it is based on the Min-Max characterization of the eigenvalues of self-adjoint compact operators. \\

Concerning the computational techniques of the eigenfunctions $\varphi^{\alpha}_{n,c}(\cdot)$ and the eigenvalues $\nu_{n,\alpha}(c),$ we use a classical  efficient  and robust method based on the computation of the eigenfunctions of a   Sturm-Liouville differential operator, commuting with the integral operator $\mathcal L_c^{\alpha}.$ This technique is well known in the literature. It has been used by D. Slepian and his co-authors to compute the prolate spheroidal wave functions (PSWFs) that are concentrated on an interval in the 1-D case or on the unit disk in the 2-D case, see \cite{Slepian1, Slepian2}. These PSWFs are related to the eigenfunctions of the finite Fourier transform and the finite  Hankel transform, for the 1-D and 2-D case, respectively. Note that  various generalizations of the PSWFs can be found  in the literature, see for example \cite{ Wang2, Abreu, Moumni, Zayed}. \\

 By a  straightforward modification of the weighted finite Fourier transform case, it can be easily checked that the operator $\mathcal L_c^{\alpha}$ commutes with the following Sturm-Liouville differential operator,     
\begin{equation}\label{Sturm_operator}
\mathcal{D}_c^{\alpha}(f)(x)= \frac{d}{ d x} \left[\omega_{\alpha}(x)(1-x^2) f'(x)\right]+(c^2 x^2) \omega_{\alpha}(x) f(x).
\end{equation}
In the sequel, we let $\chi_{n,\alpha}(c)$ denote the eigenvalue of $-\mathcal{D}_c^{\alpha}$, associated with the eigenfunction $\varphi_{n,c }^{\alpha},$ so that we have
\begin{equation}\label{Sturm_operator2}
-\mathcal{D}_c^{\alpha}(\varphi_{n,c}^{\alpha})(x)=- \frac{d}{ d x} \left[\omega_{\alpha}(x)(1-x^2) (\varphi_{n,c}^{\alpha})'(x)\right]-(c^2 x^2)\omega_{\alpha}(x) \varphi_{n,c}^{\alpha}(x)=
\chi_{n,\alpha}(c) \varphi_{n,c}^{\alpha}(x),\quad x\in I.
\end{equation}
Since the differential operator is nothing but a perturbation of the Jacobi polynomials differential operator, by the quantity $c^2 x^2 \omega_{\alpha}(x),$ then we check that  a practical scheme for the computation of the $\varphi_{n,c}^{\alpha}$
is given by the series expansion of this later in the basis of orthonormal Jacobi polynomials.\\

We emphasize on the fact that the study of the decay and the behaviour of the spectrum of the finite and weighted finite Laplace transform is not explored yet in the literature. One of the aims of this work is to provide some first results concerning the fast decay of the eigenvalues of the weighted finite Laplace transform operator, as well as an estimate of the growth for the largest eigenvalue of this operator. It is important to point out that
the techniques, used in the previous joint works of one of us, \cite{Karoui-Souabni,  Karoui-Souabni2, Bonami-Karoui3}  to study the spectrum of the 
singular values of the finite and weighted finite Fourier transform operators, cannot be used in our present case. This is due to the fact that 
unlike the Fourier case, the singular values of the weighted finite Laplace transform are not related to an energy maximization problem, nor they satisfy
a simple ordinary differential equation. To overcome these limitations, we develop a Min-Max based technique to get the previous two results
concerning the spectrum of the weighted finite Laplace transform.\\

This work is organized as follows. In section~2, we give some mathematical preliminaries and describe the computational scheme of the eigenfunctions $\varphi_{n,c}^{\alpha}$ and their associated eigenvalues $\nu_{n\alpha}(c).$ In section 3, we study some qualitative and quantitative behaviours of the eigenvalues $\nu_{n,\alpha}(c).$ In particular, we prove the super-exponential decay rate of these eigenvalues. In section 4, we study some bounds of the $\varphi_{n,c}^{\alpha}$ and show that they are well adapted for the spectral approximation of a special set of functions. In section 5, we provide the reader with some numerical examples that illustrate the different results of this work.

\section{Mathematical Preliminaries and Computation of the eigenfunctions and the eigenvalues.}

In this paragraph, we give some mathematical preliminaries on some special functions that are frequently used in this paper. Also, we describe a Jacobi based method for the computation of the eigenfunctions and the eigenvalues of the weighted finite Laplace transform. We first recall that  $\J_k,$  the Jacobi polynomial of degree $k$ and parameters $(\alpha,\alpha)$ are 
normalized so that ${\displaystyle \J_k(1)=\frac{1}{k B(k,\alpha+1)}},$ where $B(\cdot,\cdot)$ is the Beta function, given by ${\displaystyle B(x,y)=\frac{\Gamma(x) \Gamma(y)}{\Gamma(x+y)}.}$ Here, $\Gamma(\cdot)$ is the Gamma function, that satisfies the following bounds, see \cite{Batir} that 
\begin{equation}\label{Eq2.2}
    \sqrt{2e}\left(\frac{x+\frac{1}{2}}{e}\right)^{x+\frac{1}{2}}\leq \Gamma(x+1) \leq  \sqrt{2\pi}\left(\frac{x+\frac{1}{2}}{e}\right)^{x+\frac{1}{2}},\quad x>-1.
\end{equation} 

Moreover,  it is known that, see for example [\cite{NIST}, p. 457] 
 \begin{eqnarray}\label{Eq2.1}
    \dint_{-1}^1 e^{cxy}\J_k(y) \omega_{\alpha}(y)d y&=&\frac{e^{-c x}2^{2\alpha+k+1}}{k!}(c x)^k B(k+\alpha+1,k+\alpha+1)\nonumber\\
   && \qquad\qquad\times     _1F_1(k+\alpha+1,2k+2\alpha+2, 2c x),\quad x\in I.
  \end{eqnarray}
where $_1F_1(\cdot,\cdot, \cdot)$ is the Kummer's  function.  It is well known that for $Re(b)> Re(a)>0$ and $z\in \mathbb C,$  this hypergeometric function
has the following integral representation,
\begin{equation}\label{Kummer}
_1F_1(a,b,z)= \frac{\Gamma(b)}{\Gamma(a) \Gamma(b-a)} \int_0^1 e^{zt} t^{a-1} (1-t)^{b-a-1}\, dt.
\end{equation}
On the other hand, the identity \eqref{Eq2.1} can be written in terms of the modified Bessel function of the first  type
$I_{\mu}, \mu>-1.$ In fact, by using \eqref{Eq2.1} and the identity
$$_1F_1(\beta+1/2,2\beta+1, 2 x)= 2^{\beta} \Gamma(\beta+1)\frac{e^{x}}{(x)^{\beta}} I_{\beta}(x),$$
one gets 
\begin{equation}\label{Eq2.4}
    \dint_{-1}^1 e^{cxy}\J_k(y) \omega_{\alpha}(y)d y = 2^{2k+3\alpha+3/2}\frac{\Gamma(k+\alpha+3/2)}{k!} B(k+\alpha+1,k+\alpha+1)\frac{I_{k+\alpha+1/2}(cx)}{(cx)^{\alpha+1/2}}
  \end{equation}
It is well known that the orthonormal Jacobi polynomials $\wJ_n,$ $n\geq 0$ given by 
\begin{equation}\label{JacobiP}
\wJ_{n}(x)= \frac{1}{\sqrt{h_n^{\alpha}}}\J_n(x),\quad h_n^{\alpha}=\frac{2^{2\alpha+1}(\Gamma(n+\alpha+1))^2}{n!(2n+2\alpha+1)\Gamma(n+2\alpha+1)}.
\end{equation}
form an orthonormal basis of the Hilbert space $E=L^2(I,\omega_{\alpha}).$ Moreover, since for any integer $n\geq 0,$ 
$\varphi_{n,c}^{\alpha}\in E,$ then we have the series expansion
\begin{equation}\label{expansion1}
    \varphi_{n,c}^{\alpha}(x)=\ds\sum_{k\geq 0} d_{k}^n \, \wJ_k(x),\quad x\in I.
\end{equation}
Since the $\wJ_k$ are the eigenfunctions of the non perturbed  Sturm-Liouville operator 
$-\mathcal D_0^{\alpha},$ with the associated eigenvalue $\chi_{k,\alpha}(0)=k(k+2\alpha+1),$ then by using a technique similar to the one used in \cite{Karoui-Souabni} and based on the three-terms recurrence relation of Jacobi polynomials,
one can easily check that the expansion coefficients $(d_k^n)_k$ as well as the associated eigenvalue $\chi_{n,\alpha}(c)$ are computed by solving the following three diagonal eigen-system.

\begin{eqnarray}\label{eigensystem}
\lefteqn{-\frac{\sqrt{(k+1)(k+2)(k+2\alpha+1)(k+2\alpha+2)}}{(2k+2\alpha+3)\sqrt{(2k+2\alpha+5)(2k+2\alpha +1)}} c^2 d_{k+2}^n }\nonumber \\
&&\qquad\qquad+\left( k(k+2\alpha+1)-c^2 \frac{2k(k+2\alpha+1)+2\alpha-1}{(2k+2\alpha+3)(2k+2\alpha-1)} \right)
d_k^n  \nonumber  \\
&&\qquad\qquad -\frac{\sqrt{k(k-1)(k+2\alpha)(k+2\alpha-1)}}{(2k+2\alpha-1)\sqrt{(2k+2\alpha+1)(2k+2\alpha-3)}} c^2
d_{k-2}^n= \chi_{n,\alpha}(c) d_k^n, \quad k\geq 0.
\end{eqnarray}
From the general spectral theory of Sturm-Liouville operators, it is easy to check that the eigenfunctions $\varphi_{n,c }^{\alpha} $ share some properties with the Jacobi polynomials, in particular, the set $\mathcal{B}=\{{\varphi_{n,c }^{\alpha} , n \geq 0}\}$ is an orthonormal basis of $L^2(I,\omega_{\alpha })$. Moreover, for any integer $n\geq 0,$
  $\varphi_{n,c }^{\alpha}$ has the same parity as $n,$ that is 
       \begin{equation}\label{17}
        \varphi_{n,c }^{\alpha}(-x)= (-1)^n\varphi_{n,c }^{\alpha }(x),  \quad \forall\; x\in I.
       \end{equation}
Also we should mention that the $(n+1)-$th eigenvalue $\chi_{n,\alpha}(c)$ satisfies the following classical inequalities
\begin{equation}\label{vp}
    n(n+2\alpha+1)-c^2 \leq \chi_{n,\alpha}(c) \leq n(n+2\alpha+1).
\end{equation}
To get the previous upper and lower  bounds, we consider the following form of the Sturm-Liouville operator, associated with the $\varphi_{n,c}^{\alpha},$
\begin{eqnarray}
-\mathcal{D}_c^{\alpha}(f)(x)&=&-\frac{1}{\omega_{\alpha}(x)}\left[\omega_{\alpha}(x)f'(x)(1-x^2)\right]'-(c^2x^2)f(x)\cr
&=&-\mathcal{D}_0^{\alpha}(f)(x)-(c^2x^2)f(x).
\end{eqnarray}
Then, from the well known  Min-Max characterization of the eigenvalues of a self-adjoint operator applied to the operator $-\mathcal{D}_c^{\alpha}$ with $G_{n+1} \subset L^2(I,\omega_{\alpha}),$ an   $(n+1)-$dimensional subspace, one gets
\begin{eqnarray*}\label{min}
\chi_{n,\alpha}(c) &=& {\ds\min_{H\in G_{n+1} }\max_{\begin{array}{c}
                                                             U\in H \\
                                                             \|U\|_{L^2(I,\omega_{\alpha})}=1
                                                           \end{array}
    }\int_{-1}^1}-\mathcal{D}_c^{\alpha}(U)(x) U(x) \omega_{\alpha}(x)d x\cr
&=& \ds{\min_{H\in G_{n+1} }\max_{\begin{array}{c}
                                                             U\in H \\
                                                             \|U\|_{L^2(I,\omega_{\alpha})}=1
                                                           \end{array}
    }}\left( \ds\int_{-1}^{1}-\mathcal{D}_0^{\alpha}(U)(x) U(x) \omega_{\alpha}(x)d x+\ds\int_{-1}^{1}(-c^2x^2)U(x) U(x) \omega_{\alpha}(x)dx\right)
\end{eqnarray*}
The inequalities given by \eqref{vp} follow from the facts that ${\displaystyle -c^2\leq \int_{-1}^{1}(-c^2x^2)U(x) U(x) \omega_{\alpha}(x)dx\leq 0}$ and the Min-Max characterization of $\chi_{n,\alpha}(0)=n(n+2\alpha+1).$ On the other hand, it is interesting to note that for a fixed integer $n\geq 0,$ the set of the  expansion coefficients $(d_k^n)_k$ has a fast decay to zero. This is given by the following proposition.

\begin{proposition}\label{decay_coefficients}
Let  $c > 0$, $\alpha >-1$ be two real numbers and let  $n, k \geq 0$ be two integers. Let 
 \begin{equation}\label{expansion_coefficient1}
d_{k}^n=\dint_{-1}^1 \wJ_k(x) \varphi_{n,c}^{\alpha}(x) \omega_{\alpha}(x)\, dx
 \end{equation}
 Then, we have
\begin{equation}\label{expansion_coefficient2}
|d_{k}^n|\leq \dfrac{C_{\alpha}}{\nu_{n,\alpha}(c)}\dfrac{e^c}{k+\frac{1}{2}}\left(\dfrac{ec}{2k+1}\right)^k,\qquad 
C_{\alpha}=\frac{2^{-\alpha}\pi^{3/2}}{\sqrt{e(1+\alpha)}}.
\end{equation}
\end{proposition}

\noindent
{\bf Proof:}
By combining \eqref{Eq1.2}, the identity \eqref{Eq2.4} and   \eqref{expansion_coefficient1},  one gets 
\begin{eqnarray*}
|d_{k}^n|&=&\dfrac{1}{\nu_{n,\alpha}(c)\sqrt{h_k^{\alpha}}}\left|\dint_{-1}^1 \J_k (x) \dint_{-1}^1 e^{cxy}\varphi_{n,c}^{\alpha}(y) \omega_{\alpha}(y) dy\omega_{\alpha}(x) dx \right|\nonumber \\
   &=& \dfrac{1}{\nu_{n,\alpha}(c)\sqrt{h_k^{\alpha}}} \left|\dint_{-1}^1 \left( \dint_{-1}^1 e^{cxy}\J_k (x) \omega_{\alpha}(x)\, dx \right)\varphi_{n,c}^{\alpha}(y) \omega_{\alpha}(y) dy \right|\nonumber \\
 &\leq & \dfrac{2^{2k+3\alpha+3/2}}{\nu_{n,\alpha}(c)\sqrt{h_k^{\alpha}}}\frac{\Gamma(k+\alpha+3/2)}{k!}B(k+\alpha+1,k+\alpha+1) \int_{-1}^1 \frac{\left|I_{k+\alpha+1/2}(cy)\right|}{|cy|^{\alpha+1/2}} |\varphi_{n,c}^{\alpha}(y)| \omega_{\alpha}(y)\, dy.
\end{eqnarray*}
It is well known that for $\mu>-1,$ see for example [\cite{NIST}, p.227]
$$ |J_{\mu}(z)|\leq \frac{|z|^{\mu}}{2^\mu \Gamma(\mu+1)} e^{|Im(z)|},\quad z\in \mathbb C.$$
On the other hand, it is well known, see [\cite{NIST}, p. 251] that the Bessel function and the modified Bessel function of the first kinds  $J_{\mu},$ and $I_{\mu}$ are related to each others by the following relation
$$I_{\mu}(z) =e^{-i \pi \mu/2} J_{\mu}(iz),\quad  -\pi\leq  \mbox{Arg} z \leq \frac{\pi}{2}.$$
By using the previous two inequalities as well as the previous  identity, one gets 
$$ |d_{k}^n|\leq \dfrac{2^{2k+3\alpha+3/2}}{\nu_{n,\alpha}(c)\sqrt{h_k^{\alpha}}}\frac{c^k}{k!}\frac{B(k+\alpha+1,k+\alpha+1)}{2^{k+\alpha+1/2}} e^c
\int_{-1}^1 |y|^k |\varphi_{n,c}^{\alpha}(y)| \omega_{\alpha}(y)\, dy.$$ Since ${\displaystyle \|\varphi_{n,c}^{\alpha}\|_{L^2(I,\omega_{\alpha})}=1},$ then by H\"older's inequality, one gets 
\begin{equation}\label{Eq2.5}
 |d_{k}^n | \leq e^c \dfrac{2^{k+2\alpha+1}}{\nu_{n,\alpha}(c)\sqrt{h_k^{\alpha}}}\frac{c^k}{k!}B(k+\alpha+1,k+\alpha+1) 
 \sqrt{B(k+1/2,\alpha+1)}. 
\end{equation}
Note that from the useful inequalities of  the Gamma function, given by  \eqref{Eq2.2}, one gets 
\begin{equation}\label{Eq2.6}
  B(k+\alpha+1,k+\alpha+1) \leq  \dfrac{\sqrt{2}\pi}{\sqrt{2k+2\alpha+3/2}2^{2k+2\alpha+1}}
\end{equation}
On the other hand,  from \cite{Alzer}, we have 
\begin{equation}\label{Eq2.7}
  B(k+1/2,\alpha+1) \leq  \dfrac{1}{(k+1/2)(\alpha+1)}.
\end{equation}
Also, by using the concavity of the logarithmic function and the inequalities \eqref{Eq2.2}, one can  check that
$$2 \frac{(\Gamma(k+\alpha+1))^2}{\Gamma(k+2\alpha+1) k!} \geq \frac{e}{\pi},$$ so that 
\begin{equation}\label{Eq2.8}
\frac{1}{\sqrt{h_k^{\alpha}}} \leq  \sqrt{\frac{\pi}{e}}\sqrt{2k+2\alpha+1}\,\, 2^{-\alpha}.
\end{equation}
Finally, since from \eqref{Eq2.2},we have ${\displaystyle \frac{c^k}{k!}\leq \frac{1}{\sqrt{2k+1}}\left(\dfrac{2ec}{2k+1}\right)^k,}$ then  by combining \eqref{Eq2.5}, \eqref{Eq2.6}, \eqref{Eq2.7} and \eqref{Eq2.8}, one gets the desired result
\eqref{expansion_coefficient2}.\qquad $\Box$\\

\noindent
Next, it is a classical result that the identity \eqref{Eq1.2} gives the analytic extension of the $\varphi_{n,c}^{\alpha}$
to the whole $\mathbb R^*_+$ and by symmetry (recall that $\varphi_{n,c}^{\alpha}$ has the same parity as $n$), this extension holds over the whole $\mathbb R^*.$ This extension as well as the explicit expression of the eigenvalue $\nu_{n,\alpha}(c)$ are given by the following lemma.

\begin{lemma} Under the notation of the previous proposition, the analytic extension of $\varphi_{n,c}^{\alpha}(x)$ is given by
\begin{equation}\label{Eq2.9}
\varphi_{n,c}^{\alpha}(x)=\frac{1}{\nu_{n,\alpha}(c)}  \sum_{k=0}^{+\infty} d_k^n \frac{2^{2k+3\alpha+3/2}\Gamma(k+\alpha+3/2)}{\sqrt{h_k^{\alpha}} k!} B(k+\alpha+1,k+\alpha+1)\frac{I_{k+\alpha+1/2}(c x)}{(cx)^{\alpha+1/2}},\quad x\in \mathbb R^*,
\end{equation}
where
\begin{equation}\label{Eq2.10}
\nu_{n,\alpha}(c)=\frac{1}{\left(\ds{\sum_{k=0}^{+\infty}}\frac{d_k^n }{k B(\alpha+1,k)\sqrt{h_k^{\alpha}}}\right) } \ds{\sum_{k=0}^{+\infty}}d_k^n \frac{2^{2k+3\alpha+3/2}\Gamma(k+\alpha+3/2)}{\sqrt{h_k^{\alpha}} k!} B(k+\alpha+1,k+\alpha+1)\frac{I_{k+\alpha+1/2}(c )}{c^{\alpha+1/2}}.
\end{equation}
\end{lemma}

\noindent
{\bf Proof:} To prove \eqref{Eq2.9}, it suffices to use \eqref{Eq1.2} and insert the expansion \eqref{expansion1} in the integral, given in the identity \eqref{Eq2.4}. Then by using the fast decay of the expansion coefficients $(d_k^n)_k,$ given by the previous proposition, together with the fact that, see for example [\cite{NIST}, p.450]
$$\max_{x\in [-1,1]} |\J_k(x)|=\J_k(1)=\frac{1}{k B(k,\alpha+1)},\quad \alpha>-1,$$
one can interchange the integral and the sum signs and obtain the identity \eqref{Eq2.9}. Finally, the explicit expression of 
$\nu_{n,\alpha}(c),$ given by \eqref{Eq2.10} follows from equalling at $x=1,$ the two  expansions of $\varphi_{n,c}^{\alpha},$ given by \eqref{expansion1} and \eqref{Eq2.9}.$\qquad \Box$\\

\section{ Behaviour and  fast decay rate of the eigenvalues of the weighted 
finite Laplace transform.}

In this paragraph, we use the Min-Max theorem and show that the sequence of the eigenvalues $\nu_{n,\alpha}(c)$ has a super-exponential decay rate to zero. This is given by the following theorem, which is one of the main results of this 
work.

\begin{theorem}
 For given real numbers $c > 0$,$ \;\; \;\alpha>-1$ and for any integer $n > \frac{ec}{2}+1,$ we have
\begin{equation}\label{Eq3.1}
    \nu_{n,\alpha}(c)\leq e^{c} \frac{k_{\alpha}}{\log\left(\frac{2n-1}{ec}\right)} \exp\left(-(n-1)\log\left(\frac{2n-1}{ec}\right)\right),\quad k_{\alpha}=\frac{2^{-\alpha}\pi^{3/2}}{\sqrt{ e(\alpha+1)}}.
\end{equation}
\end{theorem}

\noindent
{\bf Proof:} We first recall the Courant-Fischer-Weyl Min-Max variational principle concerning the positive eigenvalues
of a self-adjoint compact operator $A$ on a Hilbert space $\mathcal H,$ with positive eigenvalues arranged in the decreasing order $\lambda_0\geq \lambda_1\geq \cdots\geq \lambda_n\geq \cdots ,$ then we have 
$$\lambda_n = \min_{f\in S_n}\,\,\,  \max_{f\in S_n^{\perp},\|f\|_{\mathcal H}=1} < Af, f>_{\mathcal H},$$
 where $S_n$ is a subspace of $\mathcal H$ of dimension $n.$ In our case, we have $A=\mathcal L_c^{\alpha},$ $\mathcal H= L^2(I,\omega_{\alpha}).$ We consider the special case of  $$S_n=\mbox{Span}\left\{\wJ_0, \wJ_1,\ldots,\wJ_{n-1}\right\}$$ and
 $$f=\ds\sum_{k\geq n} a_k \wJ_k \in S_n^{\perp},\qquad 
  \parallel f\parallel_{L^2_{(I,\omega_\alpha)}}=\ds\sum_{k\geq n}|a_k|^2=1.$$
From the proof of proposition~\ref{decay_coefficients}, we have 
\begin{equation}\label{Eq3.2}
\| \mathcal L_c^{\alpha} \wJ_k\|_{L^2(I,\omega_{\alpha})} \leq k_{\alpha} e^c \left(\frac{ec}{2k+1}\right)^{k},\quad k_{\alpha} =\frac{2^{-\alpha}\pi^{3/2}}{\sqrt{e(\alpha+1)}}.
\end{equation}
Hence, for the previous $f\in S_n^{\perp},$  and by using H\"older's inequality, combined with  the Minkowski's inequality for an infinite sum and taking into account that  $\|f\|_{L^2(I,\omega_{\alpha})}= 1,$ so that $|a_k|\leq 1,$ for $k\geq n,$  one gets
\begin{equation}\label{Eq3.3}
|< \mathcal L_c^{\alpha} f,f>_{L^2(I,\omega_{\alpha})}| \leq \sum_{k\geq n} |a_k| \| \mathcal L_c^{\alpha} \wJ_k\|_{L^2(I,\omega_{\alpha})}\leq e^c k_{\alpha} \sum_{k\geq n} \left(\frac{ec}{2k+1}\right)^{k}.
\end{equation}
The decay of the sequence appearing in the previous sum, allows us to compare this later with its integral counterpart, that is 
\begin{equation}\label{Eq3.4}
\sum_{k\geq n} \left(\frac{ec}{2k+1}\right)^{k}\leq \int_{n-1}^{\infty} e^{-x \log(\frac{2x+1}{ec})}\, dx\leq 
\int_{n-1}^{\infty} e^{-x \log(\frac{2n-1}{ec})}\, dx.
\end{equation}
Hence, by using \eqref{Eq3.3} and \eqref{Eq3.4}, one concludes that 
\begin{equation}\label{Eq3.5}
\max_{f\in S_n^{\perp},\, \|f\|_{L^2(I,\omega_{\alpha})=1}} < \mathcal L_c^{\alpha} f,f>_{L^2(I,\omega_{\alpha})}\leq e^{c} \frac{k_{\alpha}}{\log(\frac{2n-1}{ec})} e^{-(n-1)\log(\frac{2n-1}{ec})},\quad k_{\alpha}=\frac{2^{-\alpha}\pi^{3/2}}{\sqrt{ e(\alpha+1)}}.
\end{equation}
 To conclude for the proof of the theorem, it suffices to use the previous  Courant-Fischer-Weyl Min-Max variational principle.$\qquad \Box $
 
 \begin{remark}
 A result similar to the result of the previous theorem, but for the case of the singular values of the finite Fourier transform operator,  has been recently given  in \cite{Bonami-Karoui5}. Moreover, we should mention that in this last special case, a much more elaborated highly accurate and explicit approximation formula, leading to a sharp  super-exponential decay rate of the singular values has been given in \cite{Bonami-Karoui3}. 
  \end{remark}
 
It is interesting to note that even for large values of $c>0,$  and as shown by the previous theorem, the sequence of the eigenvalues $\nu_{n,\alpha}(c)$ decays super-exponentially to zero, the first and largest eigenvalue $\nu_{0,\alpha}(c)$  has large value  which is comparable with $e^c.$ This behaviour of the large value is given by the following proposition.

\begin{proposition}
For real numbers $c>0,$ $\alpha > -\frac{1}{2}$ and any $0<\gamma<1,$ we have 
\begin{equation}\label{Eq3.6}
\nu_{0,\alpha}(c) \geq K_{\alpha,\gamma} e^{(\gamma^2 c)},\qquad K_{\alpha,\gamma}= \sqrt{\frac{2}{\pi}}\frac{(1-\gamma)^2(1-\gamma^2)^{\alpha-1/2}}{(\alpha+1)}.
\end{equation}
\end{proposition}

\noindent
{\bf Proof:} By the Max-Min Theorem, we have 
$$\nu_{0,\alpha}(c)=\max_{S_1}\min_{f\in S_1} <\mathcal L_c^{\alpha} f, f>_{L^2(I,\omega_{\alpha})},$$
where $S_1$ is a subspace of $L^2(I,\omega_{\alpha})$ of dimension $1.$ In particular for $S_1=\mbox{Span} \{ \wJ_0\},$ we have
$$\nu_{0,\alpha}(c)\geq \,\,\,\, \big <\mathcal L_c^{\alpha} \wJ_0, \wJ_0\big >_{L^2(I,\omega_{\alpha})}.$$
Taking into account that 
$$\mathcal L_c^{\alpha} \wJ_0(x)=\frac{2^{3\alpha+3/2}}{\sqrt{h_0^{\alpha}}}\Gamma(\alpha+3/2) B(\alpha+1,\alpha+1)\frac{I_{\alpha+1/2}(cx)}{(cx)^{\alpha+1/2}},\qquad h_0^{\alpha} =2^{2\alpha+1} B(\alpha+1,\alpha+1),$$
one gets 
$$ \mathcal L_c^{\alpha} \wJ_0(x) = 2^{2\alpha+1} \Gamma(\alpha+3/2)\sqrt{B(\alpha+1,\alpha+1)}\frac{I_{\alpha+1/2}(cx)}{(cx)^{\alpha+1/2}},\quad x>0.
$$
Moreover, since ${\displaystyle \wJ_0(x)= \frac{1}{\sqrt{2^{2\alpha+1} B(\alpha+1,\alpha+1)}}}$ and by using the parity of $\mathcal L_c^{\alpha} \wJ_0(x),$ one gets 
\begin{equation}\label{Eq3.7}
\big <\mathcal L_c^{\alpha} \wJ_0, \wJ_0\big >_{L^2(I,\omega_{\alpha})}=2^{\alpha+1/2} \Gamma(\alpha+3/2)\int_0^1 
\frac{I_{\alpha+1/2}(cx)}{(cx)^{\alpha+1/2}} \omega_{\alpha}(x)\, dx.
\end{equation}
Also, from the following integral representation of the $I_{\mu}(x),$ see for example [\cite{NIST}, p.252],
$$I_{\mu}(x)= \left(\frac{x}{2}\right)^{\mu} \frac{1}{\sqrt{\pi} \Gamma(\alpha+1/2)}\int_{-1}^1 (1-t^2)^{\mu-1/2} e^{xt} \, dt,\quad x>0,$$
one gets for $0<\gamma <1,$
$$\int_{-1}^1 (1-t^2)^{\mu-1/2} e^{xt} \, dt \geq e^{\gamma x} (1+\gamma)^{\mu-1/2}\int_{\gamma}^1(1-t)^{\mu-1/2}\, dt.$$ 
Hence, by combining the previous inequality and equality, one gets 
\begin{equation}\label{Eq3.8}
I_{\alpha}(x) \geq C_{\alpha,\gamma} x^{\alpha} e^{\gamma x},\qquad C_{\alpha,\gamma}=\frac{(1-\gamma^2)^{\alpha-1/2}(1-\gamma)}{2^{\alpha}\sqrt{\pi} \Gamma(\alpha+3/2)},\quad x>0.
\end{equation} 
By using the previous inequality, it is easy to check that 
\begin{equation}\label{Eq3.9}
\int_0^1 
\frac{I_{\alpha+1/2}(cx)}{(cx)^{\alpha+1/2}} \omega_{\alpha}(x)\, dx\geq C_{\alpha,\gamma} e^{\gamma^2 c}\int_{\gamma}^1 (1-x^2)^{\alpha}\, dx\geq C_{\alpha,\gamma} \frac{(1-\gamma^2)^{\alpha}}{1+\alpha}(1-\gamma) e^{(\gamma^2 c)}.
\end{equation}
Finally, by combining \eqref{Eq3.7} and \eqref{Eq3.9}, one gets the desired result \eqref{Eq3.6}.$\qquad\Box$\\

We should mention that it is easy to compute the trace of the operator  $\mathcal L_c^{\alpha}.$ In fact, 
since $ k(x,y)=e^{cxy},$ the kernel of $\mathcal L_c^{\alpha}$  is  continuous symmetric and non-negative definite, then by  Mercer's theorem, the trace of  $\mathcal{L}_c^{\alpha}$ is given by
 \begin{eqnarray*}
    Trace(\mathcal{L}_c^{\alpha})&=&\ds{\sum_{n=0}^{+\infty}}\nu_{n\alpha}(c)=\ds{\int_{-1}^{1}}e^{cx^2}(1-x^2)^{\alpha}dx.
 \end{eqnarray*}
  Also, by use the substitution $t=x^2$ and the integral representation of the Kummer's function \eqref{Kummer}, one gets 
 \begin{equation}\label{Trace}
  Trace(\mathcal{L}_c^{\alpha})=\ds{\int_{0}^{1}}e^{ct}(1-t)^{\alpha}\frac{1}{\sqrt{t}}dt=B\left(\frac{1}{2}, \alpha+1\right) \,\,   _1F_1\left(\frac{1}{2},\alpha+\frac{3}{2},c\right).
 \end{equation}

\section{Eigenfunctions of $\mathcal L_c^{\alpha}:$ Uniform bounds and their quality of approximation.}

In this paragraph, we first  give an upper bound for the eigenfunction $\varphi_{n,c}^{\alpha}.$ This bound generalizes or  improves the bounds obtained in \cite{Bonami-Karoui1, Bonami-Karoui2} and  \cite{Karoui-Souabni}, in the cases of the finite Fourier transform operator and the  weighted finite Fourier transform operator, respectively. Then, by using this bound together with the fast decay rate of the eigenvalues $\nu_{n,\alpha}(c),$ we check that the $\varphi_{n,c}^{\alpha}$ are well adapted for the approximation of functions from the space $\mathcal L_c^{\alpha} \left(L^2(I,\omega_{\alpha})\right).$ Note that to get an uniform bound of $\varphi_{n,c}^{\alpha},$ we first need a local estimate of $\varphi_{n,c}^{\alpha},$ given by the following lemma. Note that the proof of this lemma is similar, and slightly different from the proof of proposition 2 of \cite{Karoui-Souabni}, given in the weighted finite Fourier transform operator case.  

\begin{lemma}
 Let  $c > 0$,$ \;\alpha > -1,$ be two real numbers. Then, for any integer $n\in \mathbb N,$ satisfying
 ${\displaystyle q=\frac{c^2}{\chi_{n,\alpha}(c)}\geq 0}$ and $q\leq \frac{\alpha}{2}$ if $\alpha>0,$ we have:
 \begin{equation}\label{local_estimate}
 \ds \sup_{t\in[-1,1]}(1-t^2) \omega_{\alpha}(t)\big(\varphi_{n,c}^{\alpha}(t)\big)^2
   \leq 1+\alpha.
 \end{equation}
\end{lemma}

\noindent
{\bf Proof: }  The proof uses a classical technique for the local estimates of the eigenfunctions of a Sturm-Liouville operator. In our case, the eigenfunction $\varphi_{n,c}^{\alpha}$ has the same parity as $n.$
Hence, it suffices to prove the previous estimate on the interval $[0,1].$ We  consider the auxiliary function, defined on $[0,1] $ by
\begin{equation}\label{Auxiliary}
     Z_n(t)=\big(\varphi_{n,c}^{\alpha}(t)\big)^2+\frac{1-t^2}{(1+
     qt^2){\chi_{n,\alpha}(c)}}\big((\varphi_{n,c}^{\alpha})'(t)\big)^2.
\end{equation}
 Since
$$
    -\mathcal{D}^{(\alpha )}_c (\varphi_{n,c}^{\alpha}(t))=\chi_{n,\alpha}(c) \varphi_{n,c}^{\alpha}(t),
$$
then a straightforward computation gives us,
 \begin{equation}\label{20}
     Z'_n(t)=\frac{2((\varphi_{n,c}^{\alpha})')^2(t)}{\chi_{n,\alpha}(c)(1+ qt^2)}
     \left[ t(2\alpha+1)-\frac{2t(1-t^2)q}{(1+ qt^2)} \right].
\end{equation}
Next, we consider a second auxiliary function, given by
$$
    K_n(t)=(1-t^2)\omega_{\alpha}(t)Z_n(t),\hspace{4cm}\forall\; t\in [0,1].
$$
By using the previous expression of $Z'_n(t),$ one can easily check that 
 $$
    K_n'(t)= -2t(\alpha+1)\omega_{\alpha}(t)(\varphi_{n,c}^{\alpha}(t))^2+ A(t)((\varphi_{n,c}^{\alpha})'(t))^2\frac{\omega_{\alpha}(t)(1-t^2)}{\chi_{n,\alpha}(c)(1+ qt^2)},
$$
where 
\begin{equation}
\label{Eqq4.1}
A(t)= 2t \left[\alpha -\frac{2q(1-t^2)}{(1+qt^2)}\right]= \frac{2t}{(1+qt^2)}\left[(\alpha+2)q t^2+(\alpha-2q)\right].
\end{equation}
We first consider the case  $\alpha > 0,$ then by using the previous second form of $A(t),$ one concludes that 
$A(t)\geq 0,$ for any $t\in [0,1]$ and any $n \in \mathbb{N}$ with  $0\leq q \leq \frac{\alpha}{2}.$ Hence, in this case, we have
$$
   K_n'(t) \geq  -2t(\alpha+1)\omega_{\alpha}(t)(\varphi_{n,c}^{\alpha}(t))^2,\quad t\in [0,1].
$$
Moreover,  since $K_n(1) =0$ and since $\ds 2 \int_{0}^1(\varphi_{nc}^{\alpha}(t))^2\omega_{\alpha}(t)dt=1$,
then from the previous inequality, one gets 
$$
\sup_{x\in [0,1]} (1-x^2) \omega_{\alpha}(x) \big(\varphi_{n,c}^{\alpha}(x)\big)^2\leq 
 \sup_{x\in [0,1]} K_n(x)\leq 2 (\alpha+1)  \int_{0}^1  \omega_{\alpha}(t)(\varphi_{n,c}^{\alpha}(t))^2 dt=\alpha+1.
$$
Finally, if $-1<\alpha<0,$ then by considering the interval $[-1,0]$ and using the first form of the quantity 
$A(t),$ given by \eqref{Eqq4.1}, one concludes that $A(t)\geq 0,$ for $t\in [-1,0].$ Again, since $K_n(-1)=0,$ then by applying the same steps as in the case where $\alpha >0,$ one concludes that 
$$\sup_{x\in [-1,0]} (1-x^2) \omega_{\alpha}(x) \big(\varphi_{n,c}^{\alpha}(x)\big)^2\leq \alpha+1,$$
whenever $-1<\alpha <0.$ This concludes the proof of the lemma.$\qquad \Box$

Once the local estimate \eqref{local_estimate} has been established, we prove the following theorem that provides us 
with a uniform  bound of the eigenfunction $\varphi_{n,c}^{\alpha}.$ 

\begin{theorem}
Let  $c > 0$,$ \;\; \;\alpha> -\frac{1}{2}$ be two real numbers. Then, for any integer $n\in \mathbb N,$ satisfying
 $${\displaystyle \chi_{n,\alpha}(c) \geq 6 \left(\frac{\alpha+1}{\alpha+3}\right),\qquad \mbox{ and } q=\frac{c^2}{\chi_{n,\alpha}(c)} \leq \left\{ \begin{array}{ll} \frac{1}{2}+\alpha
&\mbox{ if } -1/2<\alpha \leq 0\\
\min\left(1,\frac{1}{2}+\alpha\right)
&\mbox{ if } \alpha >0,\end{array}\right.}$$
  we have 
   \begin{equation}\label{Bound1}
   \ds\sup_{t\in [-1,1]}|\varphi_{n,c}^{\alpha}(t)|=
   |\varphi_{n,c}^{\alpha}(1)|\leq 3 \sqrt{\alpha+1} \left(\chi_{n,\alpha}(c)\right)^{\frac{\alpha+1}{2}}.
    \end{equation}
\end{theorem}

\noindent
{\bf Proof:} We first note that since $\varphi_{n,c}^{\alpha}$ has the same parity as $n,$ then it suffices to
have a bound over $[0,1]$. For this purpose, it is a classical technique to consider the same auxiliary function $Z_n(\cdot)$ given by \eqref{Auxiliary}.  Straightforward computation gives us
 \begin{equation}\label{Eq4.2}
     Z'_n(t)=\frac{2t \left((\varphi_{n,c}^{\alpha})'\right)^2(t)}{\chi_{n,\alpha}(c)(1+ qt^2)^2}
    \left[ (2\alpha+3)q t^2+(2\alpha+1-2q) \right].
\end{equation}
Hence, for any integer $n\in \mathbb N,$ with  $ 0\leq q= \frac{c^2}{\chi_{n,\alpha}(c)} \leq \alpha+\frac{1}{2},$ 
the function $Z_n(\cdot)$ is positive and  increasing. Consequently, we have   
 ${\displaystyle  \ds\sup_{[ 0,1] } Z_n(t)=Z_n(1)=(\varphi_{n,c}^{\alpha}(1))^2.}$ That is 
      \begin{equation}\label{Eq4.3}
   \ds\sup_{[ 0,1] }|\varphi_{n,c}^{\alpha}(t)|=|\varphi_{n,c}^{\alpha}(1)|
   \end{equation}
   To  bound the quantity $|\varphi_{n,c}^{\alpha}(1)|,$ we proceed as follows. Since for  $x\in [ 0,1]$ we have,
   \begin{eqnarray*}
     ((1-x^2) \omega_{\alpha}(x)(  \varphi_{n,c}^{\alpha})'(x))'&=&-\chi_{n,\alpha}(c)
     \omega_{\alpha}(x)(1+qx^2)\varphi_{n,c}^{\alpha}(x).
   \end{eqnarray*}
Then by an integration over $[x,1],$ where  $x\in [0,1[ $, one gets
\begin{eqnarray}\label{31}
    |(\varphi_{n,c}^{\alpha})'(x)|&\leq  &\frac{\chi_{n,\alpha}(c)}{(1-x^2)\omega_{\alpha}(x)}
    \ds\int_x^1 \omega_{\alpha}(t)(1+ qt^2)| \varphi_{n,c}^{\alpha}(t)| dt
    \leq \frac{\chi_{n,\alpha}(c)}{(1+x)}(1+q)| \varphi_{n,c}^{\alpha}(1)|.
\end{eqnarray}
A second integration over $[x,1],$ gives us 
\begin{eqnarray}\label{Eq4.4}
 |\varphi_{n,c}^{\alpha}(x)| &\geq& |\varphi_{n,c}^{\alpha}(1)| \left(1-\chi_{n,\alpha}(c)
  (1+ q)\log\Big(\frac{2}{1+x}\Big)\right). 
\end{eqnarray}
Let $x_n\in ]0,1[ $ be such that $${\displaystyle \log\Big(\frac{2}{1+x_n}\Big)=\frac{\gamma}{\chi_{n,\alpha}(c)},}$$ where the constant $\gamma$ is to be fixed later on. Note that from this choice of $x_n,$ we have
$$1-x_n^2= 4 \big(1-e^{-\gamma/\chi_{n,\alpha}(c)}\big) e^{-\gamma/\chi_{n,\alpha}(c)}.$$
Moreover, since
$$(1-e^{-x}) e^{-x} \geq \frac{3}{4} x,\quad 0\leq x\leq \frac{1}{6},$$
then, we have 
\begin{equation}\label{Eq4.5}
1-x_n^2 \geq 3 \frac{\gamma}{\chi_{n,\alpha}(c)},\quad\mbox{whenever \ \ } \frac{\gamma}{\chi_{n,\alpha}(c)}\leq \frac{1}{6}.
\end{equation}
By substituting $x$ with $x_n$ in \eqref{Eq4.4} and using \eqref{Eq4.5} together with the local estimate given by 
\eqref{local_estimate}, one gets 
\begin{eqnarray}
  |\varphi_{n,c}^{\alpha}(1)| &\leq & \sqrt{\alpha+1}\frac{1}{1-\gamma (1+q)}(1-x_n^2)^{-(\alpha+1)/2}\cr
  & \leq &3^{-(\alpha+1)/2} \sqrt{\alpha+1} \frac{\gamma^{-(\alpha+1)/2}}{1-\gamma (1+q)}\left(\chi_{n,\alpha}(c)\right)^{\frac{\alpha+1}{2}}.
\end{eqnarray}
The function ${\displaystyle \gamma \rightarrow \frac{\gamma^{-(\alpha+1)/2}}{1-\gamma (1+q)}}$ has the unique critical point ${\displaystyle \gamma=\frac{\alpha+1}{(3+\alpha)(1+q)}}$, which correspond to its maximum. Also, from the  conditions on $\chi_{n,\alpha}(c),$ the inequality \eqref{Eq4.5} is always satisfied.  Moreover, by substituting this value of $\gamma$ in the previous inequality and using the fact that $q\leq 1, $ one gets
\begin{equation}
  |\varphi_{n,c}^{\alpha}(1)| \leq  \sqrt{\alpha+1} \frac{(\alpha+3)}{2} \left(\frac{2(\alpha+3)}{3(\alpha+1)}\right)^{(\alpha+1)/2}\left(\chi_{n,\alpha}(c)\right)^{\frac{\alpha+1}{2}}.
\end{equation}
Finally, it is easy to check that the function defined for $\alpha \geq -1/2,$ by ${\displaystyle f(\alpha)=\frac{(\alpha+3)}{2} \left(\frac{2(\alpha+3)}{3(\alpha+1)}\right)^{(\alpha+1)/2}}$ has a global maximum at $\alpha=3$ and its maximum is given by $f(3)=3.$ This concludes the proof of the theorem. $\qquad \Box$\\

\noindent
In the last part of his paragraph, we check via the following proposition,  that the $\varphi_{n,c}^{\alpha}$ are well adapted 
for the approximation of functions from the space $\mathcal L_c^{\alpha} \left(L^2(I,\omega_{\alpha})\right).$ 
Note that the analysis of similar spectral approximations by the eigenfunctions of the finite Fourier transform operator or the more general weighted finite Fourier transform operator, have been given in \cite{Boyd, Wang1} and \cite{Karoui-Souabni, Wang2}, respectively. 

\begin{proposition}
Let  $c > 0$,$ \;\; \alpha > -1/2$ be two real numbers and let $g=\mathcal{L}_c^{(\alpha)} f,$  where $f\in L^2(I,\omega_{\alpha})$. Let $S_n(g)(x)=\ds \sum_{k\leq n} <g,\varphi_{k,c}^{\alpha}>_{L^2(I,\omega_{\alpha})} \varphi_{k,c}^{\alpha}.$ For any integer $n\geq \frac{ ec }{2}$ and satisfying the conditions of 
theorem~2, we have 
 \begin{equation}\label{Approx1}
   \|g-S_n(g)\|_{L^2(I,\omega_{\alpha})}\leq k_{\alpha}  e^c \frac{1}{\log\left(\frac{2n+1}{ec}\right)} \exp\Big(-n \log\big(\frac{2n+1}{ec}\big)\Big) \|f\|_{L^2(I,\omega_{\alpha})},\quad k_{\alpha}= \frac{2^{-\alpha}\pi^{3/2}}{\sqrt{e(\alpha+1)}}
\end{equation}
and
\begin{equation}\label{Approx2}
   \sup_{x\in [-1,1]} |g(x)-S_n(g)(x)|\leq K'_{\alpha}   e^c \frac{((n+1)(n+2\alpha+2))^{\frac{\alpha+1}{2}}}{\log\left(\frac{2n+1}{ec}\right)} \exp\Big(-n \log\big(\frac{2n+1}{ec}\big)\Big)  \|f\|_{L^2(I,\omega_{\alpha})},
\end{equation}
where  $ K'_{\alpha}$ is a constant depending only  on $\alpha$.
\end{proposition}

\noindent
{\bf Proof: } We first recall that since $\mathcal{L}_c^{\alpha}$ is a self-adjoint compact and one to one operator,
then $\mathcal{B}=\{{\varphi_{n,c }^{\alpha} , n \geq 0}\}$ is an orthonormal basis of $L^2(I,\omega_{\alpha })$. Hence, for $f\in L^2(I,\omega_{\alpha}),$ we have 
 $f(y)=\ds \sum_{n\geq 0} b_n(f) \varphi_{nc}^{\alpha}(y)$ with $b_n(f)=<f,\varphi_{nc}^{\alpha}>_{L^2(I,\omega_{\alpha})}.$ 
 Moreover, it is easy to check that if $f\in L^2(I,\omega_{\alpha}),$ then $g =\mathcal L_c^{\alpha} f \in L^2(I,\omega_{\alpha}).$ Hence, we have 
 \begin{equation}
 \label{Expansion1}
 g(x)= \sum_{k=0}^{\infty} b_k(g) \varphi_{k,c}^{\alpha}(x),\quad b_{k}(g)= \int_{-1}^1 g(t) \varphi_{k,c}^{\alpha}(t)\omega_{\alpha}(t)\, dt,\quad x\in [-1,1].
 \end{equation}
Here, the equality holds in the $L^2(I,\omega_{\alpha})-$sense. Moreover, we check  that for such a function $g,$ the corresponding previous series expansion converges absolutely to $g.$  In fact, since the function $(x,y)\rightarrow e^{cxy} f(y) \varphi_{k,c}^{\alpha}(x)$ belongs to $L^1(I\times I, \omega_\alpha(y)\omega_\alpha(x)),$ (Recall that $I$ is a compact interval), then by Fubini's theorem, we have 
\begin{eqnarray}\label{coeff}
b_k(g) &=& \int_I\left(\int_I e^{cxy} f(y)\omega_{\alpha}(y)\, dy\right) \varphi_{k,c}^{\alpha}(x)\omega_{\alpha}(x)\, dx \nonumber\\
&=& \int_I\left(\int_I e^{cxy} \varphi_{k,c}^{\alpha}(x)\omega_{\alpha}(x)\, dx\right) f(y)\, dy = \nu_{k,\alpha}(c) b_k(f).
\end{eqnarray}
Also, from Parseval's identity, we have ${\displaystyle \|f\|^2_{L^2(I,\omega_{\alpha})}=\sum_{k=0}^{\infty} |b_k(f)|^2.}$ Hence, the identity \eqref{coeff} implies that
\begin{equation}
|b_k(g)| \leq \nu_{k,\alpha}(c) \| f\|_{L^2(I,\omega_{\alpha})},\quad k\geq 0.
\end{equation}
That is, the expansion coefficients of $g$ have a super-exponential decay rate to zero, similar to the decay rate of
eigenvalues $\nu_{n,\alpha}(c).$ Moreover, we have
\begin{eqnarray*}
\|g-S_n(g)\|^2_{L^2(I,\omega_{\alpha})}&=& \sum_{k\geq n+1} |b_k(g)|^2 = \sum_{k\geq n+1} \left(\nu_{k,\alpha}(c)\right)^2 
|b_k(f)|^2 \\
&\leq & \left(\nu_{n+1,\alpha}(c)\right)^2 \sum_{k\geq n+1} |b_k(f)|^2 \leq \left(\nu_{n+1,\alpha}(c)\right)^2 \|f\|^2_{L^2(I,\omega_{\alpha})}
\end{eqnarray*}
To conclude for the proof of \eqref{Approx1}, it suffices to use the decay rate of $\nu_{n,\alpha}(c),$ given by \eqref{Eq3.1}. To prove \eqref{Approx2}, it suffices to use the fact that the expansion \eqref{Expansion1} holds point-wise. Hence, by using the bound of $\varphi_{n,c}^{\alpha},$ given by \eqref{Bound1} as well the upper bound of the 
eigenvalues $\chi_{n,\alpha}(c),$ given by \eqref{vp}, one gets for any $x\in I,$
\begin{eqnarray*}
| g(x)-S_n(g)(x)| &\leq& \sum_{k\geq n+1} |b_k(g)| \sup_{x\in I} |\varphi_{k,c}^{\alpha}(x)| \\
&\leq & 3 \sqrt{\alpha+1}\left(\sum_{k\geq n+1} \nu_{k,\alpha}(c) (\chi_{k,\alpha}(c))^{(\alpha+1)/2}\right) \|f\|_{L^2(I,\omega_{\alpha})}\\
&\leq &3 \sqrt{\alpha+1} \left(\sum_{k\geq n+1} \nu_{k,\alpha}(c) (k(k+2\alpha+1))^{(\alpha+1)/2}\right) \|f\|_{L^2(I,\omega_{\alpha})}
\end{eqnarray*}
Finally, to conclude the proof of \eqref{Approx2}, it suffices to combine the previous inequality with the super-exponential decay rate of the $\nu_{k,\alpha}(c),$ given by \eqref{Eq3.1}.$\qquad \Box$

\section{Numerical results}

In this paragraph,  we give various  numerical examples that illustrate the different  results of this work.\\

\noindent
{\bf Example 1:} In this example, we first check that a truncated version of formula \eqref{Eq2.10} is  highly accurate for computing  approximate  values $\widetilde \nu_{n,\alpha}(c)$ of the corresponding eigenvalues $\nu_{n,\alpha}(c).$
To do so, we  have considered the values of $c=\pi$ and $\alpha=-\frac{3}{4}$ and  used formula \eqref{Eq2.10} truncated to the order $N=90$ to compute the approximate values of the first $80$ eigenvalues. Then, we have 
compared the exact value of the trace of $Trace(\mathcal L_c^{\alpha})=B\left(\frac{1}{2}, \alpha+1\right) \,\,   _1F_1\left(\frac{1}{2},\alpha+\frac{3}{2},c\right),$ given by \eqref{Trace} with the approximate trace $\widetilde{Tr}(\mathcal L_c^{\alpha})=
\sum_{n=0}^{80} \widetilde \nu_{n,\alpha}(c).$ We have found  that 
$$ | Trace(\mathcal L_c^{\alpha})-\widetilde{Tr}(\mathcal L_c^{\alpha})|=2.42 E-62.$$

Next, we illustrate the results of theorem 1, concerning the super-exponential decay rate 
of the eigenvalues $\nu_{n,\alpha}(c).$ For this purpose, we have used two values of $\alpha=-\frac{3}{4}$ and $\alpha=1$
and the five values of $c= k \pi, \,\, 1\leq k\leq 5.$ Then, we have used the previous  truncated version of formula \eqref{Eq2.10}  and computed highly accurate approximate  values of $\log(\widetilde \nu_{n,\alpha}(c)),$  corresponding  to the exact values $\log(\nu_{n,\alpha}(c)),$ for different values of $n.$ The obtained numerical results are given by Figure 1.
\begin{figure}[]\hspace*{0.05cm}
 {\includegraphics[width=15.0cm,height=5.5cm]{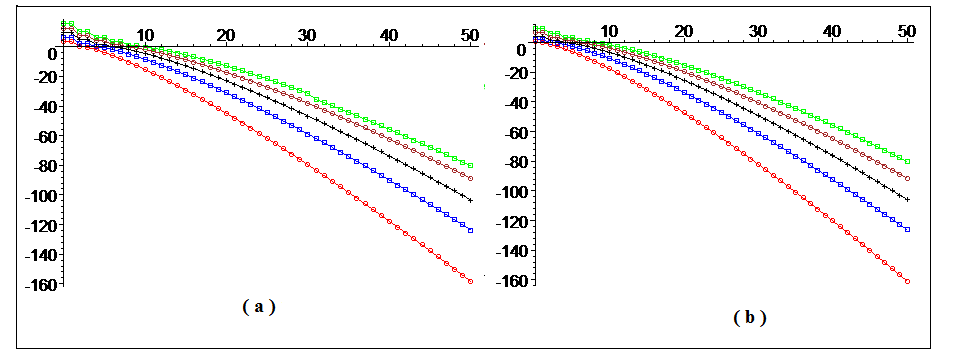}}
  \vskip -0.5cm\hspace*{1cm} \caption{  (a) Graphs of  $\log(\nu_{n,\alpha}(c))$ with $\alpha=-\frac{3}{4}$
    and different  values of $c=\pi, 2\pi, 3\pi, 4\pi, 5\pi,$ (from the left to the right), (b) same as (a) with $\alpha=1$ instead of $\alpha=-\frac{3}{4}.$   }
  \end{figure}

 Also, to check the theoretical upper bound $\overline{\nu_{n,\alpha}}(c)$ of the eigenvalue  $\nu_{n,\alpha}(c),$  given by the right hand side of \eqref{Eq3.1},  we have plotted in Figure 2, the graphs of the $\log(\nu_{n,\alpha}(c))$ versus the graphs of $\log(\overline{\nu_{n,\alpha}}(c))$ in the case where $c=\pi$ and $\alpha=-\frac{3}{4},\, 1.$ \\

 \begin{figure}[]\hspace*{0.05cm}
  {\includegraphics[width=15.0cm,height=5.5cm]{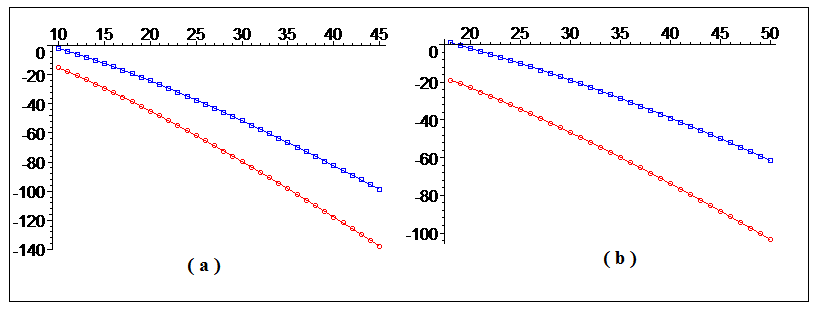}}
   \vskip -0.5cm\hspace*{1cm} \caption{  (a) Graphs of  $\log(\nu_{n,\alpha}(c))$ (in red circles) versus its theoretical bound 
   $\log(\overline{\nu_{n,\alpha}}(c))$ (in blue boxes) for $c=\pi$ and $\alpha=-\frac{3}{4},$
   (b) same as (a) with $\alpha=1$ instead of $\alpha=-\frac{3}{4}.$   }
   \end{figure}

 Next, we illustrate the property that the first eigenvalue $\nu_{0,\alpha}(c)$ has an exponential blow-up with respect to $c.$ This property has been explained by proposition 2 in the case where $\alpha > -\frac{1}{2}.$  For this purpose, we have considered the values of $c=k\pi,\,\, 1\leq k\leq 5$ and the two values of $\alpha= -\frac{3}{4},\, 1.$ Then, we have computed the values of the corresponding eigenvalues $\nu_{0,\alpha}(c),$ given by Table~1. Note that from the numerical results, this blow-up  holds for the first few values of 
$\nu_{n,\alpha}(c).$ 

\begin{center}
  \begin{table}[]
  \vskip 0.2cm\hspace*{3cm}
  \begin{tabular}{ccc} \hline
   $c$ &$\nu_{0,\alpha}(c),\, \alpha=-\frac{3}{4}$& $\nu_{0,\alpha}(c),\, \alpha=1$\\   \hline
  $\pi$  &$3.24362E+01 $   &$1.73873E+00 $  \\
  $2\pi$ &$6.19658E+02 $   &$7.51136E+00 $   \\
  $3\pi$ &$1.29094E+04 $   &$7.40701E+01 $   \\
  $4\pi$ &$2.77508E+05 $   &$9.48287E+02 $   \\
  $5\pi$ &$6.06695E+06 $   &$1.39132E+04 $    \\\hline
     & &  \\
  \end{tabular}
  \caption{Values of the first eigenvalues of $\mathcal L_c^{\alpha}$ for $\alpha=-\frac{3}{4}, 1$ and different values od $c.$}
  \end{table}
  \end{center}

 \noindent
 {\bf Example 2:} In this last example, we illustrate the spectral approximation quality by the eigenfunctions of the weighted finite Laplace transform operator.  Also, we illustrate the finite Laplace spectral based method for the inversion of this last operator.  Recall that from proposition 3, these eigenfunctions are particularly well adapted for the approximation of functions from $\mathcal L_c^{\alpha}(E),\, E=L^2(I,\omega_{\alpha}).$ Nonetheless, from the general property of  self-adjoint compact and one to one operators, the set of these eigenfunctions is also an orthonormal basis of $E.$  In this example, consider the value of $\alpha=0$ and the  test function 
 \begin{equation}\label{testfunction1}
 g_{a,\beta}(x)=\mathcal L_c^{0}(f_{a,\beta})(x),\qquad f_{a,\beta}(t) = e^{\beta t}\sin(at),\quad x\in I,
 \end{equation}
where $a,\beta\in \mathbb R.$   Straightforward computations give us the following explicit expression of $g,$
\begin{equation}\label{testfunction2}
 g_{a,\beta}(x)=\frac{2}{a^2+\beta^2+c^2 x^2+2cx\beta}\Big( (c x+\beta)\sin(a)\cosh(cx+\beta)-a\cos(a)\sinh(cx+\beta)\Big).
 \end{equation}
By using the Jacobi series expansion of the $\varphi_{n,c}^{\alpha},$ with the exponential decay rate of the expansion coefficients $(d_k^n)_k,$ combined with the identity $b_n(g)= \nu_{n,\alpha}(c) b_n(f),$ given by \eqref{coeff}, these different series expansion coefficients of $g_{a,\beta}$ are given by the following formula,
$$ b_n(g_{a,\beta}) = \nu_{n,\alpha}(c) \sum_{k\geq 0} d_k^n  <\widetilde P_k^{(0,0)}, f_{a,\beta}>_{L^2(I)}.$$

 In the special $c= a= 5 \pi$ and $\beta=3,$ we have computed $S_n(g_{a,\beta}),$ with $n=16$ and found the following uniform  and  $L^2-$errors,
given by 
$$\sup_{x\in I} \left| g_{5\pi,3}(x)- S_{16}(g_{5\pi,3})(x)\right| = 3.95E-04,\qquad \left\| g_{5\pi,3}(x)- S_{16}(g_{5\pi,3})(x)\right\|_{L^2(I)}=8.79E-05 .$$
 As predicted  by Proposition 3, our proposed  spectral approximation method based on the eigenfunctions $\varphi_{n,c}^{\alpha}$ is highly accurate for the approximation of functions from the space $\mathcal L_c^{\alpha}(L^2(I,\omega_{\alpha})).$ 
Also, it is known from the literature, see \cite{Wang2, Boyd, Wang1},   that a spectral approximation scheme for band-limited functions and based on the eigenfunctions of the finite or more generally the weighted  finite transform operators, outperforms the schemes based on the Jacobi polynomials. It seems that this is also true for the spectral approximation by the eigenfunctions of the Laplace transform. Indeed, we have repeated the previous test, by 
substituting the projection $S_{16}(g_{5\pi,3})$ by $\Pi_{16}(g_{5\pi,3}),$ the projection over the subspace spanned by the first $17$ polynomials
$\wJ_k,$ $\alpha=0.$ In this case, we have found that    
$$\sup_{x\in I} \left| g_{5\pi,3}(x)- \Pi_{16}(g_{5\pi,3})(x)\right| = 399.81,\qquad \left\| g_{5\pi,3}(x)- \Pi_{16}(g_{5\pi,3})(x)\right\|_{L^2(I)}=65.19 .$$
 We should mention that  the projection over a larger subspace, given by  $\Pi_{28}(g_{5\pi,3}),$ and spanned
by the first $29$ Jacobi polynomials   provides similar results, than the projection $S_{16}(g_{5\pi,3}).$ More precisely, we have found that  
$$\sup_{x\in I} \left| g_{5\pi,3}(x)- \Pi_{28}(g_{5\pi,3})(x)\right| = 1.82E-04,\qquad \left\| g_{5\pi,3}(x)- \Pi_{28}(g_{5\pi,3})(x)\right\|_{L^2(I)}=2.63
E-05.$$

 \noindent
 {\bf Example 3:}  In this last example, we describe the use of  the identity \eqref{coeff} for the numerical approximation of the inverse of the  weighted finite Laplace transform. In fact, if $g=\mathcal L_c^{\alpha}(f),$ where $f\in L^2(I,\omega_{\alpha})$ is unknown but for some $N\in \mathbb N,$ the projection of $g,$ given by  
${\displaystyle S_{N}(g)(x)=\sum_{k=0}^N b_k(g) \varphi_{k,c}^{\alpha}(x), \, x\in I}$ is known, then from  \eqref{coeff}, we have 
$$ f_N = S_N^{-1}(g)= \sum_{k=0}^N \frac{b_k(g)}{\nu_{k,\alpha}(c)} \varphi_{k,c}^{\alpha},$$
a numerical approximation of $f.$ Moreover, under the assumption that for some $s>0,$ the function $f$ belongs also to a Sobolev type space $\widetilde H^s,$ of those functions in $L^2(I,\omega_{\alpha})$ with a finite norm, given by 
$$\| f\|_{\widetilde H^s}^2 =\sum_{n=0}^\infty  (1+n^2)^s  |b_n(f)|^2,$$ then from the Parseval's identity, satisfied by the orthonormal basis given by the $\varphi_{k,c}^{\alpha},$ one can easily check that 
 $$\|f-f_N\|_{L^2(I,\omega_{\alpha})}^2= \|f-S_N^{-1}(g)\|_{L^2(I,\omega_{\alpha})}^2=\sum_{n\geq N} |b_n(f)|^2 \leq (N+1)^{-2s} \| f\|_{\widetilde H^s}^2.$$
 We have applied the previous scheme with $c=5\pi,$ $\alpha=0$ and  the function $g_{5\pi,3}$ of the previous example. In this case,  the exact inverse finite Laplace transform of $g_{5\pi,3}$ is  given by   $f_{5\pi,3}(x) = e^{3x}\sin(5\pi x).$ For the value of $N=30,$ we have found that 
 $$\sup_{x\in I}|f_{5\pi,3}(x)-S_N^{-1}(g_{5\pi,3})(x)|=1.67 E-04,\qquad \|f_{5\pi,3}-S_N^{-1}(g_{5\pi,3})\|_{L^2(I)} = 2.95E-05.$$

 \begin{remark}
 We should mention that due to the fast decay of the eigenvalues $\nu_{n,\alpha}(c),$ the previous scheme is numerically unstable in the presence of an added small perturbation to the function $g.$ In this case, the previous scheme has to be combined with a regularization scheme for the ill-posed problems. The study of this issue as well the issue of the performance of the spectral approximation scheme based on the eigenfunctions and eigenvalues of the weighted  finite  Laplace transform, compared to the classical scheme based on classical orthogonal polynomials, is beyond the scope of this work. These issues will be studied in a future work.
 \end{remark}


\end{document}